\documentclass[11pt]{article}
\PassOptionsToPackage{hyphens}{url}

\usepackage[margin=1in]{geometry}
\usepackage{amsmath,amssymb,amsthm}
\usepackage{graphicx}
\usepackage{booktabs}
\usepackage{array}
\usepackage{longtable}
\usepackage{multicol}
\usepackage{tikz}
\usepackage[hidelinks]{hyperref}
\usepackage{url}

\newtheorem{theorem}{Theorem}[section]
\newtheorem{lemma}[theorem]{Lemma}
\newtheorem{proposition}[theorem]{Proposition}
\newtheorem{corollary}[theorem]{Corollary}
\newtheorem{conjecture}[theorem]{Conjecture}
\theoremstyle{definition}
\newtheorem{definition}[theorem]{Definition}
\newtheorem{remark}[theorem]{Remark}
\newtheorem{example}[theorem]{Example}

\newcommand{\Z}{\mathbb{Z}}
\newcommand{\F}{\mathbb{F}}
\newcommand{\R}{\mathbb{R}}
\newcommand{\Q}{\mathbb{Q}}
\newcommand{\Hom}{\operatorname{Hom}}
\renewcommand{\Im}{\operatorname{Im}}
\newcommand{\lk}{\operatorname{lk}}
\newcommand{\tr}{\triangleright}
\newcommand{\tri}{\triangleright^{-1}}
\newcommand{\PhiE}{\Phi_E}
\newcommand{\gen}[1]{\langle #1 \rangle}

\title{Causality Detection via Symplectic Quandles}
\author{Amirbek Baxshilloyev\\
Stanford University, Stanford, CA, USA\\
\texttt{nodirovichamirbek@gmail.com}}
\date{September 2026}

\begin{document}
\maketitle

\begin{abstract}
We study whether symplectic quandle colorings can reveal causal structure encoded by
``sky links''---links formed by the spheres of all light rays through two events in the
space of light rays of a $(2+1)$-dimensional globally hyperbolic spacetime. The testbed is the
connected sum of two Hopf links $H$ (causally unrelated events) and the infinite family of
Allen--Swenberg links $L_n$, which the Alexander--Conway polynomial cannot tell apart from $H$.
We use the enhanced quandle counting polynomial, in which every coloring is weighted by the size
of the subquandle generated by its colors. For the $16$-element symplectic quandle $T=(\Z/4)^2$ we
prove that $|\Hom(Q(L_n),T)| = 640 + 96\cdot 16^{n}$ for every $n\ge 1$, while
$|\Hom(Q(H),T)|=736$, and we determine the enhanced polynomial of every $L_n$. Hence already the plain quandle
counting invariant distinguishes $H$ from every Allen--Swenberg link, and the Allen--Swenberg links
from each other. Analogous closed formulas hold for $(\Z/2)^4$ and for an $8$-element degenerate
symplectic quandle $(\Z/2)^3$ ($176+48\cdot16^n$ versus $224$), and the quandle $(\Z_7)^2$ also separates
$H$ from the whole family. Over $(\Z_p)^2$ the enhanced polynomial carries at most one number beyond
the counting invariant, and whether it separates $L_1$ from $H$ depends on $p$: it does for
$p=7,11,17,19,23,31$ and does not for $p=2,3,5,13$. The proofs combine an exact transfer-matrix
description of the family with certified finite computations.
\end{abstract}



\section{Introduction}\label{sec:intro}

Knots, in mathematical terms, are smooth embeddings of the circle $S^1$ into three-dimensional
Euclidean space $\R^3$ (or into $S^3$). Two knots are considered equivalent if they are ambient
isotopic; by Reidemeister's theorem this happens exactly when their diagrams are related by a
finite sequence of Reidemeister moves and planar isotopies. The most fundamental knot is the unknot
(or trivial knot), which is simply a geometrically unknotted loop -- a plain circle in space. When
several circles are involved, we refer to the configuration as a link; a knot is a link with one
component. Among links, the trivial link (or unlink) consists of two or more disjoint, unlinked
loops. The simplest nontrivial link is the Hopf link, formed by two unknotted circles that cross each
other twice in a diagram.

Knot invariants are functions that assign the same value to equivalent knots or links, although
distinct knots may share the same value. Such invariants have emerged as tools for understanding
causal structure in globally hyperbolic spacetimes. The Low conjecture~\cite{Low88} asserts that two
events in a $(2+1)$-dimensional globally hyperbolic spacetime whose Cauchy surface has universal cover
$\R^2$ are causally related if and only if their skies are linked in the space of light rays.
Nat\'ario and Tod~\cite{NT04} formulated the Legendrian Low conjecture in dimension $3+1$: for a
Cauchy surface homeomorphic to $\R^3$, two events are causally related if and only if their skies
are Legendrian linked.

Both conjectures were proved by Chernov and Nemirovski~\cite{CN10}; the $(2+1)$-dimensional
statement holds whenever the Cauchy surface $\Sigma$ is not homeomorphic to $S^2$ or $\R P^2$. Here
and throughout, two events are \emph{causally related} if there is a causal curve from one to the
other, that is, a piecewise smooth curve whose tangent vectors are timelike or null; informally,
information can travel between them without exceeding the speed of light. These results show a deep
connection between the topology of sky links and the causal structure of spacetime, and motivate the
use of link invariants as detectors of causality. Nat\'ario and Tod~\cite{NT04} examined specific
pairs of causally related events and observed that the associated sky links have non-trivial
Kauffman polynomial, but their results were limited to a particular class of skies.

Chernov, Martin and Petkova~\cite{CMP20} observed that Khovanov homology (a categorification of the
Jones polynomial) and Heegaard Floer homology (a categorification of the Alexander--Conway polynomial)
detect causality in $(2+1)$-dimensional globally hyperbolic spacetimes whose Cauchy surface is not
$S^2$ or $\R P^2$.

The polynomial invariants themselves are weaker than their categorifications. Allen and
Swenberg~\cite{AS21} conjectured that the Jones polynomial always detects causality in dimension
$2+1$, and they showed that the Alexander--Conway polynomial cannot distinguish 2-sky-like links from
$H$: they constructed an infinite family $\{L_n\}_{n\ge1}$ of ``2-sky-like'' links whose Conway polynomial
equals that of the connected sum of two Hopf links $H$---the sky link of causally unrelated events---while their Jones polynomials
differ from that of $H$. Recently Chernov, Harper and Kohli~\cite{CHK26} proved that the Links--Gould
polynomial distinguishes every $L_n$ from $H$. This limitation of the Alexander--Conway polynomial
underscores the need for other computable invariants that detect causal relationships.

A quandle is an algebraic structure whose axioms mirror the Reidemeister moves, which makes it a
natural source of knot and link invariants. Joyce~\cite{Joy82} and, independently,
Matveev~\cite{Mat82} associated to a link $L$ its fundamental quandle $Q(L)$; for a knot $K$ the
fundamental quandle determines $K$ up to replacing it by its reversed mirror image. A practical way of
producing computable invariants is to count quandle homomorphisms from $Q(L)$ to a fixed finite
target quandle $T$, i.e.\ to compute $|\Hom(Q(L),T)|$.

Nelson and collaborators refined this count by the enhanced quandle counting polynomial~\cite{Nel08,NN08},
which records the cardinality of the image subquandle $\Im(f)\subseteq T$ of each homomorphism
$f\in\Hom(Q(L),T)$. Symplectic quandles~\cite{NN08,Yet03} often contain rich, connected
subquandles, which makes them natural candidates for this enhancement.

Recent research has explored whether quandle invariants can detect causality in cases where the
Alexander--Conway polynomial fails. Leventhal~\cite{Lev22} showed that counting invariants for affine
Alexander quandles are insufficient; the Takasaki quandle is the special case $t=-1$ of an Alexander
quandle. Jain~\cite{Jai24} considered the symplectic quandle $(\Z_5)^2$. The first version of the
present paper considered $(\Z_3)^2$ and repeated the $(\Z_5)^2$ computations of~\cite{Jai24}.

\medskip
\noindent\textbf{What goes wrong, and what is true.} Both~\cite{Jai24} and version~1 of this paper
evaluated the enhanced polynomial by counting the number of distinct colors on the arcs of a diagram.
Since $Q(L)$ is generated by the arcs, the image of a homomorphism $f$ is the subquandle
$\gen{f(x_1),\dots,f(x_n)}$ \emph{generated} by the arc colors, which is in general larger than the
set of arc colors; only the former is a link invariant. We make the following contributions.
\begin{enumerate}
\item We give an explicit Reidemeister~II move on $H$ that changes the ``number of arc colors''
polynomial (Example~\ref{ex:R2}).
\item With the corrected definition, $\PhiE(L_1,T)=\PhiE(H,T)$ for $T=(\Z_3)^2$ and $T=(\Z_5)^2$
(Table~\ref{tab:v1}); the separation claimed in version~1 and in~\cite{Jai24} for $L_1$ disappears.
\item For $T=(\Z/4)^2$ with the standard symplectic form we prove, for all $n\ge1$,
\[
|\Hom(Q(L_n),T)| = 640+96\cdot16^n,\qquad |\Hom(Q(H),T)|=736,
\]
and we compute $\PhiE(L_n,T)$ for all $n$ (Theorem~\ref{thm:Z4}). Hence $T$ distinguishes $H$ and all
$L_n$ pairwise. For $(\Z/2)^4$ we obtain $736+480\cdot100^n$ (Theorem~\ref{thm:Z2^4}), and for the
$8$-element degenerate symplectic quandle $(\Z/2)^3$ with form $J\oplus0$ we obtain $176+48\cdot16^n$
versus $224$ for $H$ (Theorem~\ref{thm:T8}).
\item We prove a corrected ``transfer'' principle: for every coloring quandle admitting identity
colorings of the Allen--Swenberg block, $|\Hom(Q(L_n),T)|$ is non-decreasing in $n$
(Lemma~\ref{lem:transfer}). As a consequence $(\Z_7)^2$ separates $H$ from every $L_n$, and $(\Z_5)^2$
separates $H$ from every $L_n$ with $n\ge2$.
\item We verify that the link data are correct: from the Wirtinger presentations we reconstruct planar
diagrams of $L_1$, $L_2$, $L_3$ whose Kauffman-bracket spans are $116$, $216$, $316$, as predicted by
Allen and Swenberg (Section~\ref{ssec:data}).
\end{enumerate}

The structure of the paper is as follows. Section~\ref{sec:causality} reviews causal structure in
globally hyperbolic spacetimes. Section~\ref{sec:quandles} recalls quandles and the corrected enhanced
counting polynomial. Section~\ref{sec:symplectic} introduces symplectic quandles and proves the structural
facts we need. Section~\ref{sec:main} contains the computations and the main results, and
Section~\ref{sec:conclusion} concludes.

\section{Spacetime Causality}\label{sec:causality}

A spacetime $X$ is a time-oriented Lorentzian manifold. In suitable local coordinates at a point the
metric is the Lorentz product
\[
(x_1,\dots,x_n,t_1)\cdot(y_1,\dots,y_n,t_2)=x_1y_1+\dots+x_ny_n-t_1t_2 .
\]
A nonzero tangent vector $v$ is timelike if $v\cdot v<0$, null if $v\cdot v=0$ and causal if
$v\cdot v\le0$. Two points (events) $x,y\in X$ are causally related if there is a piecewise smooth
curve $\gamma$ from $x$ to $y$ whose tangent vectors are causal and all future-directed (or all
past-directed). Such a curve is called causal; it describes a trajectory that does not exceed the speed
of light.

A Cauchy surface $\Sigma$ is a subset of the spacetime such that every inextendible causal curve
intersects $\Sigma$ exactly once. Informally, $\Sigma$ carries complete information about the whole
evolution of the spacetime. A spacetime is globally hyperbolic if it has a Cauchy surface,
see~\cite{HE73}. These are the most studied spacetimes, since one version of the Strong Cosmic
Censorship conjecture of Penrose~\cite{Pen98} asserts that all physically reasonable spacetimes are
globally hyperbolic.

\begin{theorem}[Geroch~\cite{Ger70}, Bernal and S\'anchez~\cite{BS03,BS05}]
A globally hyperbolic spacetime $X$ is homeomorphic to $\Sigma\times\R$, where every slice
$\Sigma\times\{t\}$ is a Cauchy surface. Moreover $X$ is diffeomorphic to $\Sigma\times\R$ with the
$\R$-coordinate a smooth time function whose level sets are smooth spacelike Cauchy surfaces.
\end{theorem}

Geroch~\cite{Ger70} established the topological splitting, Bernal and S\'anchez~\cite{BS03}
showed that the Cauchy surface can be taken smooth and spacelike, and in~\cite{BS05} they obtained the
smooth splitting by a temporal function whose level sets are smooth spacelike Cauchy surfaces.

\begin{theorem}[Bernal and S\'anchez~\cite{BS07}]
A spacetime $X$ is globally hyperbolic if and only if
\begin{itemize}
\item it is causal, i.e.\ it contains no closed causal curves (no time travel), and
\item for all $x,y\in X$ the set $J^+(x)\cap J^-(y)$ is compact, where $J^+(x)$ and $J^-(y)$ denote
the causal future of $x$ and the causal past of $y$.
\end{itemize}
\end{theorem}

The second condition is sometimes described informally as the absence of naked singularities.

The space $\mathcal N$ of future-directed light rays (unparametrized null geodesics) of a globally
hyperbolic spacetime is identified with the spherical cotangent bundle $ST^*\Sigma$ of a Cauchy
surface; when $\Sigma=\R^2$ this is the solid torus $S^1\times\R^2$. The light rays through an event
$x$ form a sphere $S_x\subset\mathcal N$, the \emph{sky} of $x$. The sky is determined by where the
light emitted at $x$ is seen at a given moment, i.e.\ on a level set of a time function. Low's
conjecture relates causality to the linking of skies in $\mathcal N$.

\begin{theorem}[Low conjecture; proved by Chernov and Nemirovski~\cite{CN10}]
Let $(X,g)$ be a $(2+1)$-dimensional globally hyperbolic spacetime whose Cauchy surface is not
homeomorphic to $S^2$ or $\R P^2$ (equivalently, its universal cover is $\R^2$). Then two events
$x,y\in X$ are causally related if and only if their skies $S_x\cup S_y$ form a non-trivial link in
$\mathcal N$.
\end{theorem}

The implication ``causally unrelated $\Rightarrow$ skies unlinked'' is the easier half and was known
before~\cite{CN10}. In other words, causality corresponds to linking in the solid torus: $x$ and $y$
are causally related if and only if $S_x\cup S_y$ is not isotopic to the sky link of two causally
unrelated events.

\begin{theorem}[Legendrian Low conjecture~\cite{NT04}; proved in~\cite{CN10}]
Two events $x,y$ in a $(3+1)$-dimensional globally hyperbolic spacetime with Cauchy surface
homeomorphic to $\R^3$ are causally related if and only if their skies are Legendrian linked.
\end{theorem}

\section{Quandle Invariants of Links and Knots}\label{sec:quandles}

\begin{definition}
A quandle is a non-empty set $X$ with a binary operation $\tr$ satisfying
\begin{enumerate}
\item $x\tr x=x$ for all $x\in X$;
\item for every $y\in X$ the map $x\mapsto x\tr y$ is a bijection of $X$; that is, for all
$y,z\in X$ there is a unique $x\in X$ with $x\tr y=z$, which we denote $x=z\tri y$;
\item $(x\tr y)\tr z=(x\tr z)\tr(y\tr z)$ for all $x,y,z\in X$.
\end{enumerate}
\end{definition}

The first axiom is idempotency: every element acts trivially on itself. The second guarantees that
right multiplication can be undone: $(x\tr y)\tri y=x$. The third axiom is right
self-distributivity. Unlike groups, quandles are typically non-commutative and non-associative, so the
order and grouping of operations matter. The three axioms correspond to the three Reidemeister moves,
which makes quandles effective algebraic tools for distinguishing knots and links (see~\cite{Nel04}).

\begin{figure}[h]
\centering
\includegraphics[width=0.62\textwidth]{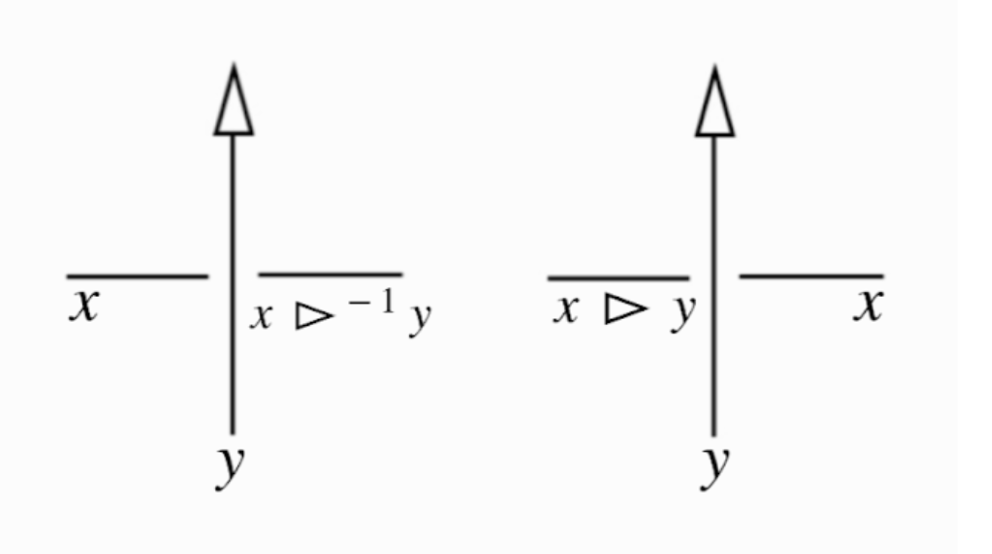}
\caption{The quandle crossing relations.}
\label{fig:crossings}
\end{figure}

Joyce~\cite{Joy82} formalized the connection between knots and quandles by assigning quandle elements
to arcs, with the operation encoding how arcs interact at each crossing, as illustrated in
Figure~\ref{fig:crossings}.

In Figure~\ref{fig:crossings}, the left diagram shows a negative crossing, where the arc $x$ passes
under $y$ from left to right, yielding $x\tri y$. The right diagram shows a positive crossing, with $x$
passing under $y$ from right to left, yielding $x\tr y$. In both cases the arc on the right of the
over-arc (with respect to its orientation) is acted on by $\tr y$ to give the arc on the left; so every
crossing can be written as a single relation $a\tr b=c$, where $b$ is the over-arc and $a$ is the
under-arc on its right. This is the form in which all relations in this paper are stored.

\begin{definition}
Let $L$ be an oriented link diagram with arc set $A$. The fundamental quandle $Q(L)$ is the quandle
generated by $A$ subject to one relation $a\tr b=c$ at every crossing.
\end{definition}

Up to isomorphism $Q(L)$ does not depend on the diagram.

\begin{definition}
Let $X_1$ and $X_2$ be quandles. A map $\varphi\colon X_1\to X_2$ is a quandle homomorphism if
$\varphi(x\tr y)=\varphi(x)\tr\varphi(y)$ for all $x,y\in X_1$. A quandle isomorphism is a bijective
homomorphism.
\end{definition}

\begin{theorem}[Joyce~\cite{Joy82}, Matveev~\cite{Mat82}]
For knots $K,K'\subset S^3$ the fundamental quandles $Q(K)$ and $Q(K')$ are isomorphic if and only if
there is a homeomorphism of $S^3$, possibly reversing orientation, taking $K$ to $K'$; equivalently,
$K'$ is equivalent to $K$ or to the reversed mirror image of $K$.
\end{theorem}

Since the fundamental quandle is hard to compare directly, one uses homomorphisms into a finite
target quandle $T$. Assigning to each arc of $L$ an element of $T$ and imposing the crossing
relations gives a system of equations whose solutions form $\Hom(Q(L),T)$.

\begin{definition}
A quandle coloring of $L$ by $T$ is an assignment of elements of $T$ to the arcs of a diagram of $L$
that satisfies $a\tr b=c$ at every crossing. Colorings are in bijection with $\Hom(Q(L),T)$, and
\[
\phi_T(L)=|\Hom(Q(L),T)|
\]
is a link invariant, the quandle counting invariant.
\end{definition}

The counting invariant ignores structural information about each homomorphism. Quandle
$2$-cocycle invariants~\cite{CEGS05} weight each homomorphism by a cocycle. Nelson~\cite{Nel08}
introduced a two-variable subquandle polynomial, and Navas and Nelson~\cite{NN08} the enhanced quandle
counting invariant, which records the size of the image of every homomorphism.

\begin{definition}\label{def:enhanced}
The enhanced quandle counting polynomial of $L$ with respect to a finite quandle $T$ is
\[
\PhiE(L,T)=\sum_{f\in\Hom(Q(L),T)}q^{|\Im(f)|}.
\]
Here $\Im(f)=f(Q(L))$ is the image of the whole fundamental quandle. It is a subquandle of $T$, and
since $Q(L)$ is generated by the arcs $x_1,\dots,x_n$ of any diagram,
\[
\Im(f)=\gen{f(x_1),\dots,f(x_n)},
\]
the subquandle of $T$ generated by the arc colors. For a finite quandle, the subquandle generated by
a set $S$ is the smallest subset containing $S$ that is closed under $\tr$ (closure under $\tri$ is
then automatic, since right multiplications are permutations of finite order).
\end{definition}

\begin{remark}\label{rem:error}
The set of colors $\{f(x_1),\dots,f(x_n)\}$ that actually appear on the arcs of a particular diagram
is in general a proper subset of $\Im(f)$, and its size depends on the diagram. Indeed, a
Reidemeister~II move creates a new arc whose color is $a\tr b$ for two existing colors $a,b$, and this
element need not have appeared before. Hence the polynomial $\sum_f q^{\#\{f(x_1),\dots,f(x_n)\}}$ is
\emph{not} a link invariant. Version~1 of this paper, and~\cite{Jai24}, computed this diagram-dependent
quantity. An explicit example is given in Example~\ref{ex:R2}. (We thank Zining Fan and
V.~Chernov for pointing this out.)
\end{remark}

Common choices of coloring quandles include:

\noindent\emph{Takasaki kei:} $\Z_n$ with $x\tr y=2y-x \pmod n$.

\noindent\emph{Affine Alexander quandle:} a module over $\Z[t,t^{-1}]$ with $x\tr y=tx+(1-t)y$. For
$t=-1$ it reduces to the Takasaki quandle.

\noindent\emph{Conjugation quandle:} a group $G$ with $x\tr y=y^{-1}xy$.

In this paper we use symplectic quandles, introduced in the next section.

\section{Symplectic Quandle}\label{sec:symplectic}

A symplectic quandle is a quandle structure on a module defined by an antisymmetric bilinear form. It
was introduced as a generalization of earlier constructions (see~\cite{NN08,Yet03}).

\begin{definition}
Let $R$ be a finite commutative ring, $M=R^r$, and let $\langle\cdot,\cdot\rangle\colon M\times M\to R$ be
a bilinear form with $\langle x,x\rangle=0$ for all $x\in M$ (so it is antisymmetric). We write
$\langle x,y\rangle=xJy^T$ for an $r\times r$ matrix $J$ with $J^T=-J$ and zero diagonal. Then $M$ is a
quandle with
\[
x\tr y=x+\langle x,y\rangle y,\qquad x\tri y=x-\langle x,y\rangle y .
\]
\end{definition}

Indeed $x\tr x=x$, and $(x\tr y)\tri y=x+\langle x,y\rangle y-\langle x+\langle x,y\rangle y,y\rangle
y=x$ because $\langle y,y\rangle=0$; self-distributivity follows from the fact that each
$R_y\colon x\mapsto x\tr y$ is a symplectic transvection, so $R_y(x\tr z)=R_y(x)\tr R_y(z)$. Unless
stated otherwise, $J$ is the standard form
$J=\begin{pmatrix}0&1\\-1&0\end{pmatrix}$ (block-diagonal for $r>2$). We write $\Z_n=\Z/n$, and
$(\Z_n)^2$ (also $(\Z/n)^2$) for the rank-$2$ symplectic quandle over $\Z_n$ with the standard form;
$(\Z/2)^4$ carries the form $J\oplus J$. Degenerate forms are allowed: the quandle
\[
T_8=(\Z/2)^3\quad\text{with}\quad \langle x,y\rangle=x_1y_2+x_2y_1\quad(\text{the form }J\oplus0)
\]
has $8$ elements and will play a role in Section~\ref{ssec:small}.

\begin{definition}
A quandle is connected if any two elements are related by a finite sequence of operations $\tr y$,
$\tri y$.
\end{definition}

\begin{theorem}[Navas and Nelson~\cite{NN08}]
Let $M$ be a symplectic quandle over a finite field with a non-degenerate form. Then the subquandle
$M\setminus\{0\}$ is connected.
\end{theorem}

The fundamental quandle of a knot is connected, while the fundamental quandle of a link with $m$
components has exactly $m$ orbits, one for each component. (Version~1 stated incorrectly that the
fundamental quandle of every nontrivial link is connected.) Over a field a non-degenerate
antisymmetric form requires even rank.

We now collect the facts about symplectic quandles that are used below.

\begin{lemma}[Scaling]\label{lem:scaling}
Let $R$ be a commutative ring and $\lambda\in R^\times$. The symplectic quandles $(R^{2k},J)$ and
$(R^{2k},\lambda J)$ are isomorphic.
\end{lemma}

\begin{proof}
Let $A=\operatorname{diag}(\lambda^{-1}I_k,I_k)$ in a symplectic basis $e_1,\dots,e_k,f_1,\dots,f_k$ with
$\langle e_i,f_j\rangle=\delta_{ij}$. Then $\langle Ax,Ay\rangle=\lambda^{-1}\langle x,y\rangle$, so
$A x+\lambda\langle Ax,Ay\rangle Ay=A(x+\langle x,y\rangle y)$, i.e.\ $x\mapsto Ax$ is an isomorphism
$(R^{2k},J)\to(R^{2k},\lambda J)$.
\end{proof}

In particular, by Lemma~\ref{lem:scaling}, the quandle $M=(\Z_3)^2$ with $J'=2J$ used in version~1 is
isomorphic to $(\Z_3)^2$ with $J$; scaling the form is not a genuine perturbation.

\begin{lemma}[Automorphisms]\label{lem:sl2}
Every $g\in SL(2,R)$ acts on $(R^2,J)$ by quandle automorphisms, and $f\mapsto g\circ f$ preserves
$\Hom(Q(L),T)$ and the size of $\Im(f)$.
\end{lemma}

\begin{proof}
$\langle gx,gy\rangle=\det(g)\langle x,y\rangle=\langle x,y\rangle$, hence $g(x\tr y)=gx\tr gy$.
\end{proof}

\begin{lemma}[The zero color]\label{lem:zero}
In any symplectic quandle, $0\tr y=0$ and $x\tr0=x$ for all $x,y$. Consequently, if one arc of a
component $K$ of $L$ is colored $0$ then every arc of $K$ is colored $0$, and colorings of $L$ with
$K\mapsto 0$ are in bijection with colorings of the sublink $L\setminus K$.
\end{lemma}

\begin{proof}
The identities are immediate from the definition. Along $K$, consecutive arcs satisfy
$x_{\rm out}=x_{\rm in}\tr^{\pm1}y$, so $x_{\rm in}=0$ forces $x_{\rm out}=0$ and conversely. At
crossings where $K$ is the over-arc the relation becomes $a\tr0=a$, which simply identifies the two
under-arcs, exactly as in a diagram of $L\setminus K$.
\end{proof}

\begin{lemma}[Subquandles of $(\F_p)^2$]\label{lem:subq}
Let $p$ be a prime and $T=(\F_p)^2$. If $S\subseteq T$ contains two linearly independent vectors, then
the subquandle generated by $S$ contains $T\setminus\{0\}$. Consequently every subquandle of $T$ is
either contained in a line $\F_p v$ (and is then a trivial quandle), or equals $T\setminus\{0\}$ or $T$.
\end{lemma}

\begin{proof}
Let $a,c\in S$ be independent, $d=\langle a,c\rangle\neq0$. The generated subquandle is invariant
under the group generated by $R_a$ and $R_c$, since $R_{x\tr y}=R_yR_xR_y^{-1}$. In the basis $(a,c)$
these are the elementary matrices $\begin{pmatrix}1&-d\\0&1\end{pmatrix}$ and
$\begin{pmatrix}1&0\\d&1\end{pmatrix}$, whose powers give all elementary matrices over the prime field
$\F_p$; they generate $SL(2,p)$, which acts transitively on $T\setminus\{0\}$. If all elements of $S$
lie on a line, then $\langle x,y\rangle=0$ for $x,y\in S$ and $S$ is already closed.
\end{proof}

In particular a subquandle of $(\Z_5)^2$ has at most $5$ elements or exactly $24$ or $25$ elements.
Image sizes such as $21$ or $22$, which appear in version~1 and in~\cite{Jai24}, cannot occur for
$\Im(f)$; they are numbers of distinct arc colors.

\section{Causality Detection using Symplectic Quandles}\label{sec:main}

\subsection{\texorpdfstring{Sky links, $H$, and the Allen--Swenberg family}{Sky links, H, and the Allen-Swenberg family}}\label{ssec:sky}

\begin{figure}[p]
\centering
\includegraphics[width=0.53\textwidth]{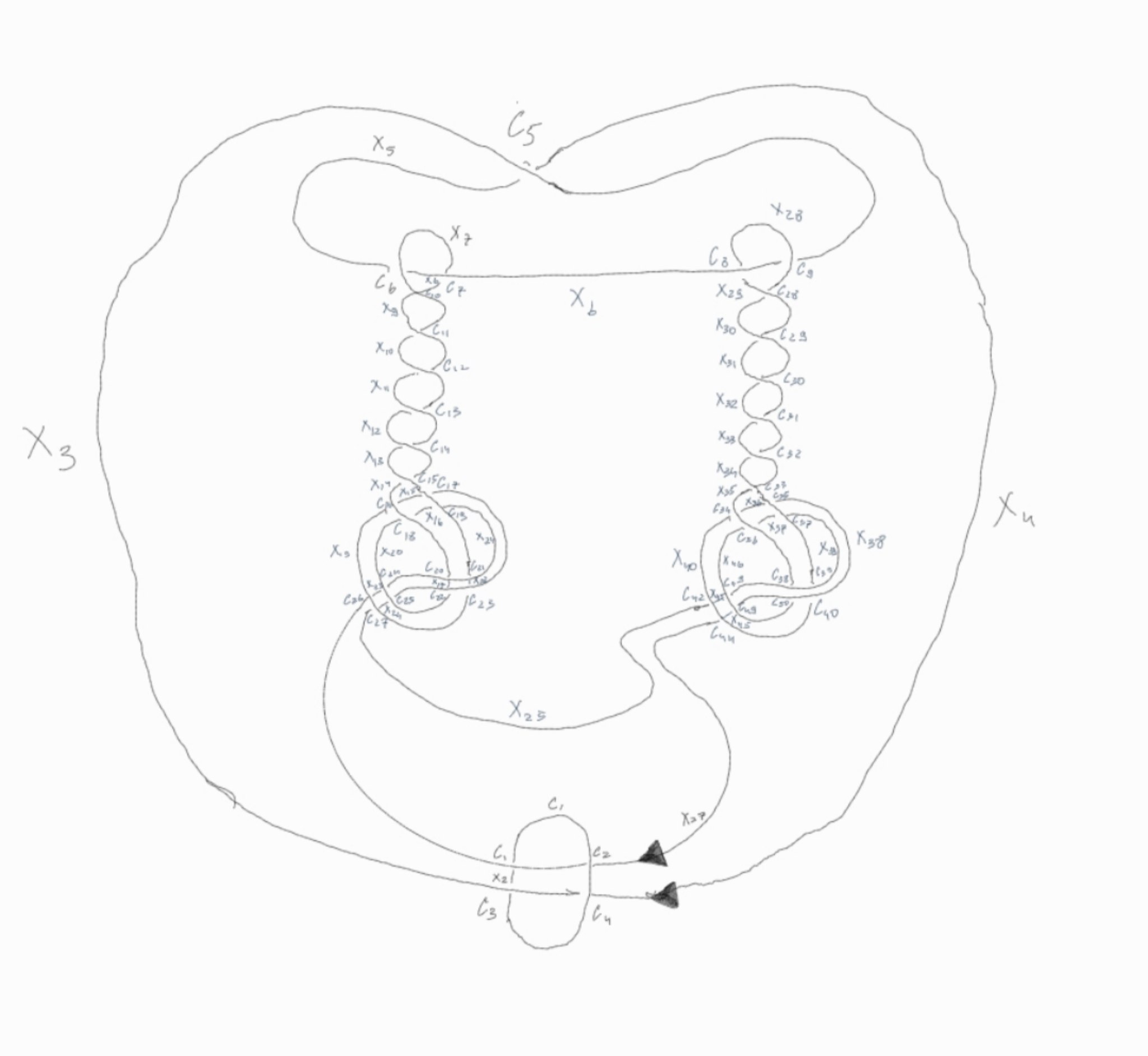}\\[2pt]
(a) Labelled Allen--Swenberg link $L_1$ (45 crossings).\\[10pt]
\includegraphics[width=0.68\textwidth]{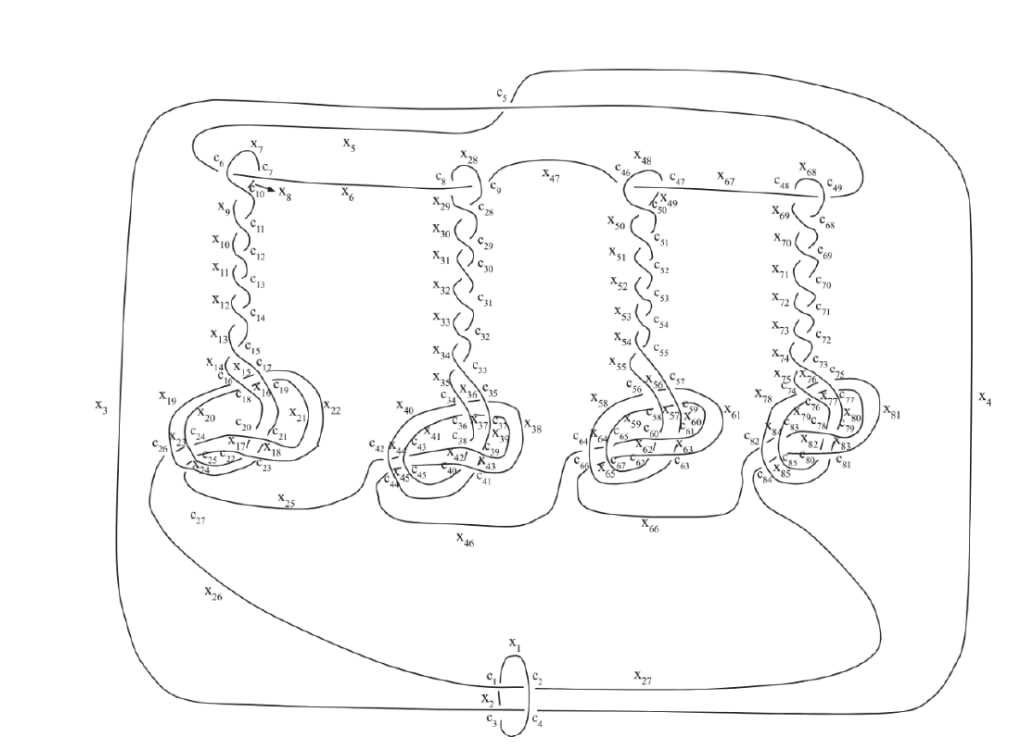}\\[2pt]
(b) Labelled Allen--Swenberg link $L_2$ (85 crossings).\\[10pt]
\includegraphics[width=0.33\textwidth]{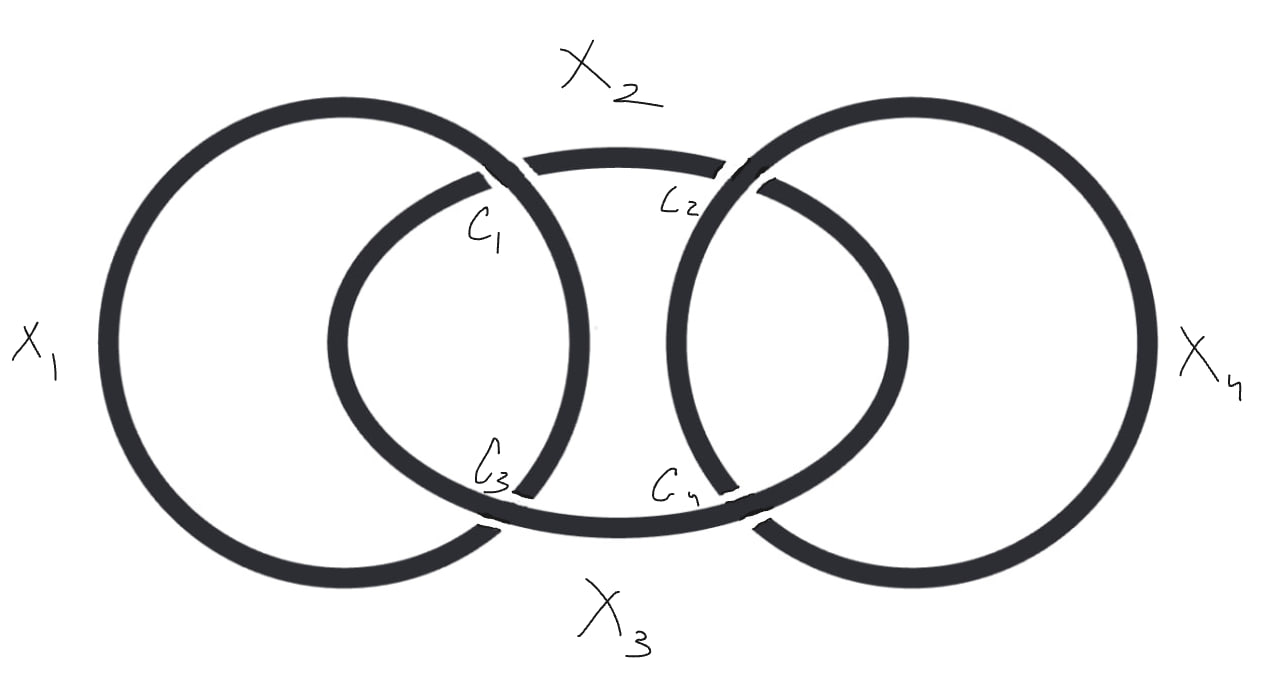}\\[2pt]
(c) Labelled connected sum of two Hopf links $H$.
\caption{The test links. In (a) and (b) the axis $\mu$ is the small circle at the bottom (arcs $x_1,x_2$).}
\label{fig:links}
\end{figure}

When the Cauchy surface is $\R^2$, the space of light rays is the solid torus $V=S^1\times\R^2$ and
the sky of an event is a circle in $V$. For two causally unrelated events the sky link is isotopic to
two parallel fibers $S^1\times\{p\}$, $S^1\times\{q\}$. To study links in $V$ with tools for links in
$S^3$ one identifies $V$ with the complement of an unknot $\mu\subset S^3$ and \emph{adds} $\mu$ as a
third component (version~1 incorrectly said ``deleting''). Each fiber links $\mu$ once and the two
fibers are unlinked, so two causally unrelated skies give the connected sum of two Hopf links
\[
H=K_1\cup\mu\cup K_2,
\]
a chain in which $\mu$ is the middle component (Figure~\ref{fig:links}(c)).

Allen and Swenberg~\cite{AS21} considered 3-component ``2-sky-like'' links $K_1\cup K_2\cup\mu$ with
$\mu$ an unknotted meridian and $\lk(K_i,\mu)=\pm1$, and constructed the infinite family
\[
L_n=C(1*T_0(n)),\qquad n\ge1,
\]
where $T_0$ is a tangle whose Conway fraction vanishes, i.e.\ $\nabla(N(T_0))=0$ and $\nabla(D(T_0))=1$ for
its numerator and denominator closures, and $T_0(n)$ is the sum of $n$ copies of $T_0$
(an explicit description of $T_0$ as a concatenation of two clasped tangles is given by Chernov,
Harper and Kohli~\cite{CHK26}). Allen and Swenberg consider two closures $C_\pm$ with the same invariants below; the checks in
Section~\ref{ssec:data} confirm that the diagrams used here (Figures~\ref{fig:links}(a),(b)) have these
invariants, and all our statements concern these diagrams.
They proved $\nabla(L_n)=\nabla(H)=z^2$ for all $n$, while the span of the Kauffman bracket of $L_n$ equals $100n+16$
and that of $H$ equals $16$. Following their proposal we treat ``causality detection''
operationally: an invariant provides evidence for detecting causality on this testbed if it separates
the Allen--Swenberg links from $H$. We work with the labelled diagrams of $L_1$ and $L_2$ in
Figures~\ref{fig:links}(a),(b) and with the family $L_n$ generated from them (Section~\ref{ssec:data}).

\subsection{Link diagrams, presentations, and their verification}\label{ssec:data}

For each oriented diagram we assign a generator $x_i$ to each arc between undercrossings and impose one
relation $a\tr b=c$ per crossing (Figure~\ref{fig:crossings}). For $H$ (Figure~\ref{fig:links}(c)):
\[
x_1\tr x_3=x_1,\quad x_3\tr x_1=x_2,\quad x_2\tr x_4=x_3,\quad x_4\tr x_3=x_4 .
\]
The presentations of $L_1$ (45 relations) and $L_2$ (85 relations) are those of the published
notebooks of version~1 (\url{https://github.com/harryfrences-svg/Horizon-Research}), read off from
Figures~\ref{fig:links}(a),(b); they are listed in Appendix~A. Because the whole paper rests on these data, we checked them independently.

\begin{enumerate}
\item \emph{Components and linking numbers.} Following under-arcs through crossings, $L_1$ has three
components with $2$ (the axis $\mu$), $39$ and $4$ arcs; $L_2$ has $2$, $77$ and $6$ arcs. Computing the
linking numbers once from the crossings where $K_i$ passes under $K_j$ and once from those where $K_j$
passes under $K_i$ gives in both cases $\lk(\mu,K_1)=\lk(\mu,K_2)=1$ and $\lk(K_1,K_2)=0$.
\item \emph{Alexander polynomials.} By Fox calculus (with $m_c=m_b^{-1}m_am_b$ for $a\tr b=c$), the
one-variable Alexander matrix of each of $H,L_1,L_2$ has a first minor equal to $(t-1)^2$ up to units,
and the multivariable Alexander polynomial of each link equals $t_\mu-1$ up to units. So the
presentations have exactly the Alexander invariants required of the Allen--Swenberg links.
\item \emph{Planar reconstruction and the Jones polynomial.} A Wirtinger presentation records the
over-arc of every crossing and which under-arc lies to its right, but not the order in which an
over-arc passes over its crossings. For $L_1$ we searched all $8\times589{,}824$ choices of component
orientations and over-crossing orders and kept those whose $4$-valent ribbon graph is planar
($V-E+F=2$). Exactly one diagram (up to reversing all orientations) is planar, so the presentation
comes from a unique classical diagram. A randomized search finds planar realizations for $L_2$ and for
the next member $L_3$ of the family (125 crossings). Their Kauffman brackets, computed by a state-sum
algorithm, have spans
\[
116=100\cdot1+16,\qquad 216=100\cdot2+16,\qquad 316=100\cdot3+16,
\]
exactly as proved by Allen and Swenberg for $L_1,L_2,L_3$. The planar diagrams are also accepted by
SnapPy~\cite{SnapPy}, which confirms that the sublinks $\mu\cup K_i$ are Hopf links.
\end{enumerate}

\emph{The family $L_n$.} The diagram of $L_2$ is that of $L_1$ with a $40$-crossing block $B$
(crossings $c_{46},\dots,c_{85}$, one copy of $T_0$) inserted between the outputs of the crossings
$c_9$, $c_{44}$ and the arcs $x_3$, $x_{27}$ that run into the clasp with $\mu$. The block is a
$4$-ended tangle with inputs $x_{47}$, $x_{46}$ and outputs $x_3$, $x_{27}$, and all its over-arcs are
internal. Inserting $n-1$ copies of $B$ produces a presentation of $L_n$; for $n=2$ this reproduces the
presentation of $L_2$ verbatim, and for $n=3$ the bracket span $316$ confirms the construction.

\begin{figure}[h]
\centering
\begin{tikzpicture}[scale=0.9,every node/.style={font=\small}]
\draw[thick] (0,0) rectangle (2.4,1.6); \node at (1.2,0.8) {$R$};
\draw[thick] (4,0) rectangle (5.4,1.6); \node at (4.7,0.8) {$B$};
\draw[thick] (6.4,0) rectangle (7.8,1.6); \node at (7.1,0.8) {$B$};
\node at (8.6,0.8) {$\cdots$};
\draw[thick] (9.4,0) rectangle (10.8,1.6); \node at (10.1,0.8) {$B$};
\draw[->] (2.4,1.2) -- (4,1.2) node[midway,above] {$x_{47}$};
\draw[->] (2.4,0.4) -- (4,0.4) node[midway,below] {$x_{46}$};
\draw (5.4,1.2) -- (6.4,1.2); \draw (5.4,0.4) -- (6.4,0.4);
\draw (7.8,1.2) -- (8.2,1.2); \draw (7.8,0.4) -- (8.2,0.4);
\draw (9.0,1.2) -- (9.4,1.2); \draw (9.0,0.4) -- (9.4,0.4);
\draw[->] (10.8,1.2) .. controls (12.2,1.2) and (12.2,2.4) .. (5.4,2.4) node[above] {$x_{3}$} .. controls (-1.2,2.4) and (-1.2,1.2) .. (0,1.2);
\draw[->] (10.8,0.4) .. controls (12.4,0.4) and (12.4,-0.8) .. (5.4,-0.8) node[below] {$x_{27}$} .. controls (-1.2,-0.8) and (-1.2,0.4) .. (0,0.4);
\end{tikzpicture}
\caption{Schematic of $L_n$: the part $R$ (crossings $c_1,\dots,c_{45}$ of $L_2$, i.e.\ $L_1$ cut open
at $x_3$, $x_{27}$) followed by $n-1$ copies of the block $B$. The picture only records how the four
boundary arcs are glued, not the geometry.}
\label{fig:transfer}
\end{figure}

\subsection{Computational method}\label{ssec:method}

All computations are exact and were done in Python~3 (package \texttt{symq}, see Reproducibility).
For a coloring quandle $T$ with operation table $\tr$ we enumerate $\Hom(Q(L),T)$ by constraint
propagation: each arc carries a set of admissible colors, each crossing $a\tr b=c$ is a constraint
that removes a color from an arc only if no triple $(x,y,x\tr y)$ with $x,y,x\tr y$ admissible uses it,
and the search branches on the arc with the fewest admissible colors. No coloring can be lost, and
every complete assignment is checked against \emph{all} crossing relations (residual test). A second
implementation of the same algorithm in C is used for speed; for $H,L_1,L_2$ over $(\Z_3)^2$, $(\Z/4)^2$
and $T_8$ both implementations return identical sets of colorings. For each coloring we compute $|\Im(f)|$ as the closure of the set of arc
colors under $\tr$ (Definition~\ref{def:enhanced}).

As a check against version~1 we also computed the diagram-dependent arc-color count: we reproduce
every number of version~1 (Table~\ref{tab:v1}),\footnote{Except the misprint $840q^3$ in the $L_1$ row of
Table~2 of version~1, whose coefficients sum to $1705\neq1225$; the text of version~1 has the correct
$360q^3$.} which shows that the difference between the two versions is due to the definition and not to
the solver.

\subsection{The error in version 1}\label{ssec:error}

\begin{example}[Reidemeister~II changes the arc-color count]\label{ex:R2}
Push the left circle of $H$ under the right circle through the outer region. The arc $x_1$ is cut into
$x_1,x_6,x_5$ and the relation $x_1\tr x_3=x_1$ is replaced by
\[
x_5\tr x_3=x_1,\qquad x_1\tr x_4=x_6,\qquad x_5\tr x_4=x_6 .
\]
Every coloring extends uniquely ($x_6=x_1\tr x_4$, $x_5=x_1$), so $|\Hom|$ and $|\Im(f)|$ are unchanged.
Over $(\Z_3)^2$, however, the 48 colorings with $\mu\mapsto0$ and linearly independent colors
$a=f(x_1)$, $c=f(x_4)$ use three colors on the original diagram and four colors on the new one, since
$a\tr c=a+\langle a,c\rangle c$ is new. The arc-color polynomial changes from $9q+72q^2+72q^3$ to
$9q+72q^2+24q^3+48q^4$, while the corrected polynomial is $9q+72q^2+24q^3+48q^9$ in both cases. Over
$(\Z_5)^2$ the arc-color polynomial changes from $25q+360q^2+840q^3$ to $25q+360q^2+360q^3+480q^4$.
\end{example}

For $H$ the corrected polynomial can be written down for every prime.

\begin{proposition}\label{prop:H}
Let $p$ be a prime and $T=(\Z_p)^2$. Then $|\Hom(Q(H),T)|=2p^4-p^2$ and
\[
\PhiE(H,T)=p^2\,q+3p(p^2-1)\,q^2+p(p-2)(p^2-1)\,q^3+p(p-1)(p^2-1)\,q^{p^2}.
\]
\end{proposition}

\begin{proof}
Over a field, $a\tr b=a$ if and only if $\langle a,b\rangle b=0$, i.e.\ $a$ and $b$ are linearly
dependent. The relations of $H$ therefore force $x_2=x_3=:b$ and say exactly that $a=f(x_1)$ and
$c=f(x_4)$ are dependent on $b$. If $b=0$ there are $p^4$ choices of $(a,c)$; if $b\neq0$ then
$a,c\in\F_pb$, giving $(p^2-1)p^2$ colorings. By Lemma~\ref{lem:subq}, $\Im(f)=\{a,b,c\}$ unless
$b=0$ and $a,c$ are independent, in which case $\Im(f)=T$ (it contains $0$); there are
$(p^2-1)(p^2-p)$ such colorings. Counting the sizes of $\{a,b,c\}$ in the remaining cases gives the
coefficients: $|\{a,b,c\}|=1$ for $p^2$ colorings, $=2$ for $3(p^2-1)+3(p-1)(p^2-1)$, and $=3$ for
$(p^2-1)(p-2)+(p^2-1)(p-1)(p-2)$.
\end{proof}

\begin{table}[h]
\centering\footnotesize
\begin{tabular}{@{}llll@{}}
\toprule
$T$ & $L$ & arc-color count (v1; not invariant) & corrected $\PhiE(L,T)$ \\
\midrule
$(\Z_3)^2$ & $H$   & $9q+72q^2+72q^3$ & $9q+72q^2+24q^3+48q^9$ \\
$|\Hom|=153$ & $L_1$ & $9q+72q^2+24q^3+48q^9$ & $9q+72q^2+24q^3+48q^9$ \\
           & $L_2$ & $9q+72q^2+24q^3+48q^9$ & $9q+72q^2+24q^3+48q^9$ \\
\midrule
$(\Z_5)^2$ & $H$ & $25q+360q^2+840q^3$ & $25q+360q^2+360q^3+480q^{25}$ \\
$|\Hom|$: & $L_1$ & $25q+360q^2+360q^3+360q^{21}+120q^{22}$ & $25q+360q^2+360q^3+480q^{25}$ \\
$1225,1225,$ & $L_2$ & $\cdots+360q^{21}+120q^{22}+12000q^{25}$ & $25q+360q^2+360q^3+12480q^{25}$ \\
$13225$ & & & \\
\bottomrule
\end{tabular}
\caption{Version~1 versus the corrected invariant. The middle column reproduces version~1 (for $L_2$
over $(\Z_5)^2$, where version~1 reported no value, the first terms are as for $L_1$). With the correct
definition $L_1$ and $H$ have the same polynomial over $(\Z_3)^2$ and $(\Z_5)^2$.}
\label{tab:v1}
\end{table}

Table~\ref{tab:v1} shows that the conclusion of version~1 does not survive: over $(\Z_3)^2$ nothing is
separated, and over $(\Z_5)^2$ the corrected invariant separates $L_2$ from $H$ and from $L_1$ (already
the plain count does so) but not $L_1$ from $H$. The same applies to the $(\Z_5)^2$ computations
of~\cite{Jai24}, which report exactly the arc-color polynomial of $L_1$ in the middle column. (For $L_2$,
\cite{Jai24} reports $|\Hom|=13125$ and a coefficient $1200q^{25}$; our computation gives
$13225=25+360+360+360+120+12000$, so these appear to be misprints.)

The mechanism is transparent. By Lemma~\ref{lem:zero}, the colorings with $\mu\mapsto0$ are the
colorings of the sky sublink $K_1\cup K_2$. For $H$ this sublink is the unlink and its colorings are the
pairs $(a,c)\in T^2$. For $L_1$ over $(\Z_3)^2$ and $(\Z_5)^2$ the number of such colorings is also
$|T|^2$; the only effect of the more complicated diagram is that the same colorings use many more
distinct colors on the arcs, which version~1 mistook for a difference of the invariant.

\subsection{\texorpdfstring{Corrected results over $(\Z_p)^2$}{Corrected results over (Z_p)^2}}\label{ssec:primes}

By Lemma~\ref{lem:subq}, over $(\Z_p)^2$ the corrected polynomial of any link can only differ from
that of $H$ in the coefficients of $q^{p^2-1}$ and $q^{p^2}$, because colorings whose image lies in a line
are colorings by a trivial quandle and their number depends only on the number of components. So the
enhancement adds at most one number to the count, namely how the non-collinear colorings split between
$q^{p^2-1}$ and $q^{p^2}$. In every case where we enumerated them ($p\le13$ for $L_1$, $p\le7$ for $L_2$),
the colorings of $L$ with $\mu\mapsto v\neq0$ are exactly those of $H$ (all colors lie on the line
$\F_pv$), so every additional coloring has $\mu\mapsto0$ and uses the color $0$. In these cases
\[
\PhiE(L,(\Z_p)^2)=\PhiE(H,(\Z_p)^2)+\bigl(|\Hom(Q(L),T)|-|\Hom(Q(H),T)|\bigr)\,q^{p^2}\qquad(L=L_1,L_2),
\]
and the corrected enhancement adds nothing to the counting invariant. Note also that the
unlink of two components has exactly $|T|^2$ colorings; since $L_1$ has more colorings with $\mu\mapsto0$
over $(\Z_7)^2$, its sky sublink $K_1\cup K_2$ is not the unlink. Table~\ref{tab:primes}
lists the counts. We count the colorings with $\mu\mapsto0$ using Lemma~\ref{lem:sl2}: fixing the arc
$x_7$ of $K_1$ to be $0$ or $e_1=(1,0)$, their number is $N_0+(p^2-1)N_{e_1}$, where $N_0$ and $N_{e_1}$ are
the numbers of colorings with $\mu\mapsto0$ and $x_7\mapsto0$, resp.\ $x_7\mapsto e_1$ (both are enumerated;
$N_0=p^2$ in all cases).

\begin{table}[h]
\centering\small
\begin{tabular}{@{}rrrrrrr@{}}
\toprule
 & \multicolumn{2}{c}{$H$} & \multicolumn{2}{c}{$L_1$} & \multicolumn{2}{c}{$L_2$}\\
\cmidrule(lr){2-3}\cmidrule(lr){4-5}\cmidrule(lr){6-7}
$p$ & $\mu\mapsto0$ & $|\Hom|$ & $\mu\mapsto0$ & $|\Hom|$ & $\mu\mapsto0$ & $|\Hom|$\\
\midrule
2 & 16 & 28 & 16 & 28 & 16 & 28 \\
3 & 81 & 153 & 81 & 153 & 81 & 153 \\
5 & 625 & 1\,225 & 625 & 1\,225 & 12\,625 & 13\,225 \\
7 & 2\,401 & 4\,753 & 120\,673 & 123\,025 & 7\,464\,289 & 7\,466\,641 \\
11 & 14\,641 & 29\,161 & 458\,161 & 472\,681 & 75\,413\,041 & $\ge$\,75\,413\,041 \\
13 & 28\,561 & 56\,953 & 28\,561 & 56\,953 & -- & -- \\
17 & 83\,521 & 166\,753 & 12\,147\,265 & $\ge$\,12\,147\,265 & -- & -- \\
19 & 130\,321 & 260\,281 & 11\,074\,321 & $\ge$\,11\,074\,321 & -- & -- \\
23 & 279\,841 & 559\,153 & 84\,219\,169 & $\ge$\,84\,219\,169 & -- & -- \\
29 & 707\,281 & 1\,413\,721 & 707\,281 & -- & -- & -- \\
31 & 923\,521 & 1\,846\,081 & 58\,062\,721 & $\ge$\,58\,062\,721 & -- & -- \\
\bottomrule
\end{tabular}
\caption{Counting invariants over $(\Z_p)^2$. The columns ``$\mu\mapsto0$'' count colorings with the axis
colored $0$, i.e.\ colorings of the sky sublink $K_1\cup K_2$ (Lemma~\ref{lem:zero}); for $H$ this is
$p^4$. Where the colorings with $\mu\mapsto v\ne0$ were not enumerated, only the lower bound
$|\Hom|\ge\#\{\mu\mapsto0\}$ is given; it already exceeds $|\Hom(Q(H),T)|$, except for $p=29$, where
the colorings with $\mu\mapsto0$ number $p^4$ as for $H$ and the rest was not enumerated. For $p\le13$ ($L_1$) and $p\le7$
($L_2$) the colorings with $\mu\mapsto v\ne0$ number exactly $(p^2-1)p^2$, as for $H$, and all their colors
lie on the line $\F_pv$.}
\label{tab:primes}
\end{table}

\begin{theorem}\label{thm:Zp}
The counting invariant of the rank-$2$ symplectic quandle $(\Z_p)^2$ distinguishes $L_1$ from $H$ for
$p=7,11,17,19,23,31$ and does not for $p=2,3,5,13$. It distinguishes $L_2$ from $H$ for $p=5,7,11$ and does
not for $p=2,3$. Over $(\Z_7)^2$ the corrected enhanced polynomials are
\begin{align*}
\PhiE(H)&=49q+1008q^2+1680q^3+2016q^{49},\\
\PhiE(L_1)&=49q+1008q^2+1680q^3+120288q^{49},\\
\PhiE(L_2)&=49q+1008q^2+1680q^3+7463904q^{49}.
\end{align*}
\end{theorem}

The set of primes that detect $L_1$ is not an interval, and it does not seem to contain all large
primes: for $p=29$ the colorings of $L_1$ with $\mu\mapsto0$ number exactly $29^4$, as for $H$ (the
colorings with $\mu\mapsto v\neq0$ were not enumerated for this $p$), while for $p=31$ the colorings with $\mu\mapsto0$ alone already outnumber all colorings of $H$. In all our data the additional
colorings are colorings of the sky sublink $K_1\cup K_2$, which is therefore not the unlink; see
Conjecture~\ref{conj:primes} for a possible explanation of the exceptional primes.

\subsection{Transfer along the Allen--Swenberg family}\label{ssec:transfer}

Version~1 claimed that every coloring of $L_1$ transfers to $L_n$ ``with the same image'' and deduced
that the enhanced polynomial separates $H$ from the whole family. As stated, the argument uses the
number of arc colors and a degree bound $\deg\PhiE(H)=3$, both of which fail for the correct invariant
(e.g.\ $\deg\PhiE(H,(\Z_5)^2)=25=\deg\PhiE(L_1,(\Z_5)^2)$). There is also a gap in the construction:
copying the colors of the first copy of $T_0$ into a new copy gives a coloring of the new block whose
input colors are those \emph{entering} the first copy, while inserting a block requires a coloring whose
inputs and outputs are both equal to the colors \emph{leaving} it. This is exactly the hypothesis of
Lemma~\ref{lem:transfer} below. The idea can be made precise as follows.

Cut $L_2$ into the part $R$ (crossings $c_1,\dots,c_{45}$) and the block $B$ (crossings
$c_{46},\dots,c_{85}$) along the four boundary arcs (Figure~\ref{fig:transfer}). For a finite quandle
$T$ and boundary colors $w=(w_1,w_2)$, $u=(u_1,u_2)\in T^2$ let
\begin{align*}
M_T[w,u]&=\#\{\text{colorings of }B\text{ with }(x_{47},x_{46})=w,\ (x_3,x_{27})=u\},\\
R_T[w,u]&=\#\{\text{colorings of }R\text{ with }(x_{47},x_{46})=w,\ (x_3,x_{27})=u\}.
\end{align*}
Since the block has no over-arcs in common with the rest of the diagram, colorings of $L_n$ are obtained
by gluing colorings of $R$ and of $n-1$ copies of $B$ along equal boundary colors, so
\begin{equation}\label{eq:transfer}
|\Hom(Q(L_n),T)|=\sum_{w,u\in T^2}R_T[w,u]\,(M_T^{\,n-1})[w,u]\qquad(n\ge1).
\end{equation}
For $n=1,2$ (and $n=3$ over $(\Z/4)^2$) formula~\eqref{eq:transfer} agrees with the direct counts in every
case we computed.

\emph{$H$ is the member $n=0$.} The part $R$ consists of the clasp $C$ (crossings $c_1,\dots,c_5$, where
$\mu$ meets the skies) followed by one copy of $T_0$ (crossings $c_6,\dots,c_{45}$). The presentation of
this copy is isomorphic to that of the block $B$ by a relabelling of arcs that sends its boundary arcs
$x_5,x_{26},x_{47},x_{46}$ to $x_{47},x_{46},x_3,x_{27}$; we found the isomorphism by computer and checked it
relation by relation. Let $A_T$ be the transfer matrix of the clasp, from the colors of $(x_3,x_{27})$ to those
of $(x_5,x_{26})$. Then
\begin{equation}\label{eq:transfer0}
|\Hom(Q(L_n),T)|=\operatorname{tr}\bigl(A_TM_T^{\,n}\bigr)\qquad(n\ge1).
\end{equation}
For $n=0$ the right-hand side counts colorings of the clasp closed up by $x_5=x_3$, $x_{26}=x_{27}$. Its
relations become $x_1\tr x_{26}=x_2$, $x_{26}\tr x_1=x_{26}$, $x_2\tr x_3=x_1$, $x_4\tr x_1=x_3$ and
$x_3\tr x_3=x_4$; the last one is a kink ($x_4=x_3$), and the others are the relations of $H$ under
$x_1\mapsto x_3$, $x_2\mapsto x_2$, $x_{26}\mapsto x_1$, $x_3\mapsto x_4$. So~\eqref{eq:transfer0} holds for
all $n\ge0$ with $L_0:=H$. (We also verified $\operatorname{tr}A_T=|\Hom(Q(H),T)|$ and
$R_T^{\mathsf T}=A_TM_T$ directly for $(\Z_3)^2$, $(\Z_5)^2$, $(\Z/4)^2$, $(\F_4)^2$, $(\Z/2)^4$ and $T_8$.)

\begin{lemma}[Corrected transfer lemma]\label{lem:transfer}
Suppose that for every $w\in T^2$ the block $B$ has a coloring with input colors $w$ and output colors
$w$, i.e.\ $M_T[w,w]\ge1$. Then $|\Hom(Q(L_{n+1}),T)|\ge|\Hom(Q(L_n),T)|$ for every $n\ge0$ (with $L_0=H$).
\end{lemma}

\begin{proof}
Given a coloring of $L_n$, let $w$ be the colors of the two arcs entering the clasp (for $n=1$ these are
the arcs $x_3$, $x_{27}$ leaving $c_9$ and $c_{44}$). Insert one more copy of $T_0$ at this place, colored by a coloring with inputs and
outputs $w$, and keep all other colors. This is a coloring of $L_{n+1}$, and different colorings of $L_n$
give different colorings of $L_{n+1}$.
\end{proof}

The hypothesis is a finite check. We verified it for $(\Z_2)^2$, $(\Z_3)^2$, $(\Z_5)^2$, $(\Z_7)^2$,
$(\Z/4)^2$, $(\F_4)^2$, $(\Z/2)^4$ and $T_8$ (for $(\Z_p)^2$ it suffices to check one pair $w$ in each orbit of
$SL(2,p)$ on $T^2$, by Lemma~\ref{lem:sl2}).

\begin{corollary}\label{cor:Z7}
The counting invariants of $(\Z_7)^2$ and of $(\F_4)^2$ distinguish $H$ from every Allen--Swenberg link
$L_n$, $n\ge1$. The counting invariant of $(\Z_5)^2$ distinguishes $H$ from every $L_n$ with $n\ge2$.
\end{corollary}

\begin{proof}
By Lemma~\ref{lem:transfer},
\[
|\Hom(Q(L_n),(\Z_7)^2)|\ge|\Hom(Q(L_1),(\Z_7)^2)|=123025>4753=|\Hom(Q(H),(\Z_7)^2)|,
\]
$|\Hom(Q(L_n),(\F_4)^2)|\ge3376>496$, and $|\Hom(Q(L_n),(\Z_5)^2)|\ge13225>1225$ for $n\ge2$. The lower
bounds only require the explicit colorings found by the solver, each of which is checked against all
relations.
\end{proof}

Table~\ref{tab:family} lists $|\Hom(Q(L_n),T)|$ for $n\le5$ computed from~\eqref{eq:transfer}; we
computed them up to $n=8$. In every column the numbers are strictly increasing in $n$ up to $n=8$, so each
of these quandles separates $L_1,\dots,L_8$ from each other (and $(\Z_5)^2$ separates $L_2,\dots,L_8$ from $H$
and from $L_1$).

\begin{table}[h]
\centering\scriptsize
\begin{tabular}{@{}lrrrrrr@{}}
\toprule
 & $(\Z_5)^2$ & $(\Z_7)^2$ & $(\Z/4)^2$ & $(\Z/2)^4$ & $(\F_4)^2$ & $T_8$\\
\midrule
$H$ & 1\,225 & 4\,753 & 736 & 1\,216 & 496 & 224 \\ \midrule
$L_1$ & 1\,225 & 123\,025 & 2\,176 & 48\,736 & 3\,376 & 944 \\
$L_2$ & 13\,225 & 7\,466\,641 & 25\,216 & 4\,800\,736 & 156\,496 & 12\,464 \\
$L_3$ & 37\,225 & 481\,705\,105 & 393\,856 & 480\,000\,736 & 7\,335\,376 & 196\,784 \\
$L_4$ & 7\,723\,225 & 32\,264\,735\,377 & 6\,292\,096 & 48\,000\,000\,736 & 357\,795\,496 & 3\,145\,904 \\
$L_5$ & 38\,371\,225 & 2\,234\,220\,779\,665 & 100\,663\,936 & 4\,800\,000\,000\,736 & 17\,413\,785\,376 & 50\,331\,824 \\
\bottomrule
\end{tabular}
\caption{$|\Hom(Q(L_n),T)|$ from the transfer formula~\eqref{eq:transfer}. Over $(\Z_2)^2$ and $(\Z_3)^2$
all values equal those of $H$ ($28$ and $153$). The columns $(\Z/4)^2$ and $(\Z/2)^4$ are
$640+96\cdot16^n$ and $736+480\cdot100^n$, and $T_8$ is $176+48\cdot16^n$ (Theorems~\ref{thm:Z4},
\ref{thm:Z2^4}, \ref{thm:T8}).}
\label{tab:family}
\end{table}

\subsection{Small symplectic quandles that separate the whole family}\label{ssec:small}

\begin{theorem}\label{thm:Z4}
Let $T=(\Z/4)^2$ with the standard form. For every $n\ge1$,
\[
|\Hom(Q(L_n),T)|=640+96\cdot16^{\,n},\qquad |\Hom(Q(H),T)|=736=640+96\cdot16^{0},
\]
and
\[
\PhiE(L_n,T)=16q+216q^2+192q^3+120q^4+48q^5+48q^6+96\cdot16^{\,n}\,q^{13},
\]
where the same formula with $n=0$ gives $\PhiE(H,T)$. In particular the plain counting invariant of
this $16$-element quandle distinguishes $H,L_1,L_2,\dots$ pairwise.
\end{theorem}

\begin{proof}
The transfer matrix $M_T$ is a $256\times256$ integer matrix. We computed $M_T$ and $R_T$ exactly and
verified, in exact integer arithmetic, that
\[
M_T^2=17M_T-16I,
\]
so the minimal polynomial of $M_T$ divides $(x-1)(x-16)$ (it equals it, since $M_T\neq I,16I$).
By~\eqref{eq:transfer0} the sequence $c_n=|\Hom(Q(L_n),T)|$ satisfies $c_{n+2}=17c_{n+1}-16c_n$ for all
$n\ge0$, and so does $640+96\cdot16^n$. Both sequences start with $c_0=736$ and $c_1=2176$ (the value for
$H$ also agrees with a brute-force count over all $16^4$ assignments), hence they coincide.

For the enhanced polynomial we enumerated the $187$ subquandles $S\subseteq T$ ($87$ trivial and
$100$ non-trivial). For a trivial $S$, $|\Hom(Q(L_n),S)|=|S|^3$. For each non-trivial $S$ we computed
the transfer matrices of colorings with values in $S$ and their exact minimal polynomials, which all
divide $(x-1)(x-16)$. Writing $e_n(S)=\#\{f:\Im(f)=S\}$, we have
$|\Hom(Q(L_n),S)|=\sum_{S'\subseteq S}e_n(S')$, and M\"obius inversion over the lattice of subquandles
gives $e_n(S)$ for all $n$ as a combination of $1^n$ and $16^n$. Summing over $|S|=m$ gives the
coefficients above, which we checked for $n\le12$; since they satisfy the same recurrence, they hold for
all $n$. Direct computations of $\PhiE$ for $H,L_1,L_2,L_3$ agree.
\end{proof}

The $13$-element subquandle carrying all additional colorings is $\{0\}\cup\{$vectors with a unit
coordinate$\}$. The values at $n=0$ are those of $H$ because of~\eqref{eq:transfer0}: $H$ is literally the ``zeroth''
Allen--Swenberg link.

\begin{theorem}\label{thm:Z2^4}
Let $T=(\Z/2)^4$ with the form $J\oplus J$. For every $n\ge1$, $|\Hom(Q(L_n),T)|=736+480\cdot100^{\,n}$,
while $|\Hom(Q(H),T)|=1216=736+480\cdot100^0$. Moreover
\begin{align*}
\PhiE(L_n,T)={}&16q+360q^2+360q^3+480q^4+720(16^n-1)q^7+360(49^n-1)q^{10}\\
&+360(49^n-2\cdot16^n+1)q^{11}+360(100^n-16^n-49^n+1)q^{15}\\
&+120(100^n+3\cdot16^n-3\cdot49^n-1)q^{16},
\end{align*}
and for $n=0$ this is $\PhiE(H,T)=16q+360q^2+360q^3+480q^4$.
\end{theorem}

\begin{proof}
As for Theorem~\ref{thm:Z4}: the transfer matrix satisfies $M^2=101M-100I$, and the subquandle lattice
($405$ subquandles) gives recurrences dividing $(x-1)(x-16)(x-49)(x-100)$. The value at $n=0$ is that of
$H$ by~\eqref{eq:transfer0}.
\end{proof}

The smallest example we found is degenerate.

\begin{theorem}\label{thm:T8}
Let $T_8=(\Z/2)^3$ with the degenerate form $J\oplus0$, i.e.\ $\langle x,y\rangle=x_1y_2+x_2y_1$. For every
$n\ge0$ (with $L_0=H$),
\[
|\Hom(Q(L_n),T_8)|=176+48\cdot16^{\,n},\qquad |\Hom(Q(H),T_8)|=224=176+48\cdot16^0,
\]
and
\[
\PhiE(L_n,T_8)=8q+96q^2+72q^3+48q^4+48(16^{\,n}-1)\,q^7,
\]
Hence already this $8$-element quandle distinguishes $H,L_1,L_2,\dots$ pairwise.
\end{theorem}

\begin{proof}
The transfer matrix is $64\times64$ and satisfies $M^2=17M-16I$; $T_8$ has $59$ subquandles and all
restricted transfer matrices have minimal polynomials dividing $(x-1)(x-16)$. The argument is then the same
as for Theorem~\ref{thm:Z4}; direct computations of $\PhiE$ for $H,L_1,L_2,L_3$ agree.
\end{proof}

Writing elements of $T_8$ as $(v,t)$ with $v\in\F_2^2$, $t\in\F_2$, the operation is
$(v,t)\tr(w,s)=(v\tr w,\;t+\langle v,w\rangle s)$. So $T_8$ is an extension of the $4$-element quandle
$(\F_2)^2$, which by itself sees nothing (Table~\ref{tab:family}), by a module of the kind used for the
cocycle invariants of~\cite{CEGS05}. The detection thus comes from a ``cocycle'' layer on top of a quandle
that is blind to the Allen--Swenberg block.

For $(\F_4)^2$ the transfer matrix has minimal polynomial
$(x-1)(x-6)(x-11)(x+9)(x-41)(x-49)(x^4+x^3+16x^2+86x+521)$; the counts are listed in
Table~\ref{tab:family}. Here, unlike in all cases we computed over $(\Z_p)^2$, additional colorings also occur
with $\mu\mapsto v\ne0$ (for $L_2$).

\subsection{Other symplectic quandles}\label{ssec:other}

Table~\ref{tab:other} summarizes the other symplectic quandles we tried on $H,L_1,L_2$. The
$(\Z_2)^2$, $(\Z_3)^2$ and $(\Z/6)^2$ quandles do not separate $H$, $L_1$, $L_2$: for $(\Z_2)^2$ and $(\Z_3)^2$ the block
$B$ is ``transparent'' (its transfer matrix is the identity), and $(\Z/6)^2\cong(\Z_2)^2\times(\Z_3)^2$.
Degenerate forms behave differently from case to case: $\lambda J$ with $\lambda$ a zero divisor of $\Z/4$ or
$\Z/6$ and $J\oplus0$ over $\Z/3$ do not separate $H$, $L_1$, $L_2$, whereas $J\oplus0$ over $\Z/2$ does
(Theorem~\ref{thm:T8}).

\begin{table}[h]
\centering\small
\begin{tabular}{@{}lrrrrl@{}}
\toprule
$T$ & $|T|$ & $|\Hom(Q(H),T)|$ & $|\Hom(Q(L_1),T)|$ & $|\Hom(Q(L_2),T)|$ & separates \\
\midrule
$(\Z/4)^2$, $J$ & 16 & 736 & 2\,176 & 25\,216 & yes \\
$(\Z/4)^2$, $2J$ (degenerate) & 16 & 1\,792 & 1\,792 & 1\,792 & no \\
$(\Z/6)^2$, $J$ & 36 & 4\,284 & 4\,284 & 4\,284 & no \\
$(\Z/6)^2$, $2J$ (degenerate) & 36 & 9\,792 & 9\,792 & 9\,792 & no \\
$(\Z/6)^2$, $3J$ (degenerate) & 36 & 20\,412 & 20\,412 & 20\,412 & no \\
$(\F_4)^2$, $J$ & 16 & 496 & 3\,376 & 156\,496 & yes \\
$(\Z/2)^4$, $J\oplus J$ & 16 & 1\,216 & 48\,736 & 4\,800\,736 & yes \\
$(\Z/2)^3$, $J\oplus0$ (degenerate) & 8 & 224 & 944 & 12\,464 & yes \\
$(\Z/3)^3$, $J\oplus0$ (degenerate) & 27 & 4\,131 & 4\,131 & 4\,131 & no \\
\bottomrule
\end{tabular}
\caption{Counting invariants for further symplectic quandles. ``Separates'' means that $H$, $L_1$, $L_2$ have
pairwise different counts.}
\label{tab:other}
\end{table}

\subsection{Consequences and outlook}\label{ssec:outlook}

The corrected picture is as follows. The enhanced polynomial over $(\Z_3)^2$ and $(\Z_5)^2$ does not
separate $L_1$ from $H$, so the main claim of version~1 was an artifact of counting arc colors. Nevertheless
symplectic quandles do detect causality on the Allen--Swenberg testbed: the $16$-element quandles
$(\Z/4)^2$ and $(\Z/2)^4$ and the $8$-element quandle $T_8$ separate $H$ and all $L_n$ pairwise by the plain
counting invariant, with closed formulas valid for all $n$; $(\Z_7)^2$ and $(\F_4)^2$ separate $H$ from all
$L_n$, and $(\F_4)^2$ also separates $L_1,\dots,L_8$ from each other.
Over $(\Z_p)^2$ the corrected enhancement adds at most one number to the count, because by
Lemma~\ref{lem:subq} the image of a coloring is either a set of collinear vectors or all of $T\setminus\{0\}$
or $T$; in our computations it adds nothing. Over rings such as $\Z/4$ the lattice of subquandles is richer, and the enhancement does record
additional structure (the $13$-element subquandle in Theorem~\ref{thm:Z4}).

\begin{conjecture}\label{conj:primes}
For every $n\ge1$ the set of primes $p$ for which $(\Z_p)^2$ distinguishes $L_n$ from $H$ is infinite, and
in fact has positive density. For $n=1$ it contains $7,11,17,19,23,31$ and does not contain
$2,3,5,13$; for $n=2$ it contains $5,7,11$ and does not contain $2,3$.
\end{conjecture}

A natural approach is to identify the colorings with $\mu\mapsto0$ with parabolic representations of
the group of the sky sublink $K_1\cup K_2$ into $SL(2,\F_p)$. If the additional colorings come from a part of
this representation variety that is defined over a number field $K$ but not over $\Q$ (for instance a
finite set of points, as for the parabolic representations of two-bridge links), then they exist over
$\F_p$ exactly for the primes with a suitable splitting behaviour in $K$. By Chebotarev's density theorem
this is a set of positive density, which need not contain all large primes; this would explain exceptional
primes such as $13$ and $29$. We have not tried to determine this variety.

\medskip
{\sloppy
\noindent\textbf{Reproducibility.} All code and data for this version are openly available at
\url{https://github.com/amirbekbakhshilloev/syplectic-quandle-amirbek-updated}. The Python package
\texttt{symq} contains: the presentations of $H$, $L_1$, $L_2$ and the generator of $L_n$; the symplectic
quandles over $\Z/n$ and $\F_{p^k}$ with axiom checks; the coloring solver (Python, with a C accelerator
implementing the same algorithm); the corrected and the arc-color polynomials; the Reidemeister~II example;
the topological checks (components, linking numbers, Alexander polynomials, planar reconstruction and
Kauffman bracket); the transfer matrices and the exact minimal-polynomial certificates. The script
\texttt{reproduce.py} regenerates every number and table in this paper; its output is stored in the
directory \texttt{results}.\par}

\section{Conclusion}\label{sec:conclusion}

This paper examined whether symplectic quandle colorings can detect causal structure encoded by sky links,
using the connected sum of two Hopf links $H$ and the Allen--Swenberg links $L_n$ as a testbed. The first
version used the number of distinct colors on the arcs of a diagram in place of the size of the image of a
coloring. This quantity is not invariant under Reidemeister moves, and with the correct definition the
separation of $H$ from $L_1$ over $(\Z_3)^2$ and $(\Z_5)^2$ disappears: the extra arc colors of the
Allen--Swenberg diagrams all lie in the subquandle generated by very few colors.

After correcting the definition and verifying the link data (Alexander polynomials, a unique planar
realization of the Wirtinger presentation of $L_1$, and Kauffman-bracket spans $116,216,316$ for
$L_1,L_2,L_3$), we obtained the following.
\begin{itemize}
\item Over $(\Z_p)^2$ the enhanced polynomial adds at most one number to the counting invariant, and
nothing in our computations; the count separates $L_1$ from $H$ for $p=7,11,17,19,23,31$ but not for
$p=2,3,5,13$.
\item A corrected transfer lemma shows that counts are non-decreasing along the family whenever the
Allen--Swenberg block admits identity colorings; hence $(\Z_7)^2$ separates $H$ from every $L_n$.
\item For the $16$-element symplectic quandle $(\Z/4)^2$ we proved $|\Hom(Q(L_n),T)|=640+96\cdot16^n$ for all
$n\ge1$ and $736$ for $H$, and computed the corrected enhanced polynomial of every $L_n$; for $(\Z/2)^4$ the
count is $736+480\cdot100^n$, and for the $8$-element degenerate quandle $T_8$ it is $176+48\cdot16^n$
versus $224$.
\end{itemize}
In particular symplectic quandles, used correctly, do detect causality on the Allen--Swenberg testbed and
distinguish the Allen--Swenberg links from each other, and a very small target quandle suffices.

\paragraph{Limitations.} Our statements concern the specific diagrams of Allen and Swenberg, whose
relation to actual skies of causally related events is itself conjectural. Quandle counting is
computationally expensive for large targets; for the family $L_n$ this is overcome by the transfer
matrix, but for general sky links the transfer structure is not available. Finally,
Conjecture~\ref{conj:primes} is supported only by computations for primes up to $31$.

\section*{Acknowledgements}

Conducted in Summer 2025 as part of the Horizon Academic Research Program (Advanced Mathematics), and
revised in 2026. This project was supervised by Professor Vladimir Chernov of Dartmouth College and
Dr.\ Ryan Maguire of MIT. I thank them both for their continual support, and I thank Zining Fan and
Vladimir Chernov for finding and explaining the error in the first version.

\appendix
\section{Wirtinger presentations}

Each entry $k\colon(a,b,c)$ stands for the relation $x_a\tr x_b=x_c$ at crossing $c_k$. For $L_2$, the
crossings $c_1,\dots,c_{45}$ agree with those of $L_1$ except $c_9\colon(6,28,47)$ and
$c_{44}\colon(45,40,46)$; the block $B$ consists of $c_{46},\dots,c_{85}$.

\noindent\textbf{$L_1$ (crossings $c_1,\dots,c_{45}$).}\\[2pt]
{\footnotesize\begin{tabular}{@{}lllll@{}}
$1\colon(1,26,2)$ & $2\colon(27,1,26)$ & $3\colon(2,3,1)$ & $4\colon(4,1,3)$ & $5\colon(5,3,4)$ \\
$6\colon(6,7,5)$ & $7\colon(7,6,8)$ & $8\colon(28,6,29)$ & $9\colon(6,28,3)$ & $10\colon(8,7,9)$ \\
$11\colon(10,9,7)$ & $12\colon(9,10,11)$ & $13\colon(12,11,10)$ & $14\colon(11,12,13)$ & $15\colon(14,13,12)$ \\
$16\colon(15,14,19)$ & $17\colon(15,13,22)$ & $18\colon(16,14,20)$ & $19\colon(16,13,21)$ & $20\colon(17,21,14)$ \\
$21\colon(18,21,13)$ & $22\colon(17,22,20)$ & $23\colon(18,22,19)$ & $24\colon(23,20,21)$ & $25\colon(24,20,22)$ \\
$26\colon(23,19,26)$ & $27\colon(24,19,25)$ & $28\colon(28,29,30)$ & $29\colon(31,30,29)$ & $30\colon(30,31,32)$ \\
$31\colon(33,32,31)$ & $32\colon(32,33,34)$ & $33\colon(35,34,33)$ & $34\colon(36,35,40)$ & $35\colon(36,34,38)$ \\
$36\colon(37,35,41)$ & $37\colon(37,34,39)$ & $38\colon(42,39,35)$ & $39\colon(43,39,34)$ & $40\colon(42,38,41)$ \\
$41\colon(43,38,40)$ & $42\colon(44,40,25)$ & $43\colon(44,41,39)$ & $44\colon(45,40,27)$ & $45\colon(45,41,38)$ \\
\end{tabular}}

\medskip
\noindent\textbf{Block $B$ of $L_2$ (crossings $c_{46},\dots,c_{85}$).}\\[2pt]
{\footnotesize\begin{tabular}{@{}lllll@{}}
$46\colon(67,48,47)$ & $47\colon(48,67,49)$ & $48\colon(68,67,69)$ & $49\colon(67,68,3)$ & $50\colon(49,48,50)$ \\
$51\colon(51,50,48)$ & $52\colon(50,51,52)$ & $53\colon(53,52,51)$ & $54\colon(52,53,54)$ & $55\colon(55,54,53)$ \\
$56\colon(56,55,58)$ & $57\colon(56,54,61)$ & $58\colon(57,55,59)$ & $59\colon(57,54,60)$ & $60\colon(62,60,55)$ \\
$61\colon(63,60,54)$ & $62\colon(62,61,59)$ & $63\colon(63,61,58)$ & $64\colon(64,58,46)$ & $65\colon(64,59,60)$ \\
$66\colon(65,58,66)$ & $67\colon(65,59,61)$ & $68\colon(68,69,70)$ & $69\colon(71,70,69)$ & $70\colon(70,71,72)$ \\
$71\colon(73,72,71)$ & $72\colon(72,73,74)$ & $73\colon(75,74,73)$ & $74\colon(76,75,78)$ & $75\colon(76,74,81)$ \\
$76\colon(77,75,79)$ & $77\colon(77,74,80)$ & $78\colon(82,80,75)$ & $79\colon(83,80,74)$ & $80\colon(82,81,79)$ \\
$81\colon(83,81,78)$ & $82\colon(84,78,66)$ & $83\colon(84,79,80)$ & $84\colon(85,78,27)$ & $85\colon(85,79,81)$ \\
\end{tabular}}

\medskip

\end{document}